\documentclass[preprint,12pt]{elsarticle}

\usepackage{amssymb}
\usepackage{amsthm}
\usepackage{epsfig}
\usepackage{times}

\usepackage{amsmath}
\usepackage{graphicx} % for graphics files
\usepackage{tikz} 
\usepackage{fancyhdr} % used to put page number in header instead of footer

\usepackage{float} % for psuedocode formatting
\usepackage{xspace}

\usepackage{mathrsfs}
\usepackage[mathcal]{euscript}
\usepackage{color}
\usepackage{array}

\usepackage[pdftex]{hyperref}

\floatstyle{ruled}
\newfloat{algorithm}{h}{lop}
\floatname{algorithm}{Algorithm}

\newcommand{\Macro}{\ensuremath{\Sigma}}
\newcommand{\Sn}{\ensuremath{S_N} }

\newcommand{\ve}[1]{\ensuremath{\mathbf{#1}}}

\newcommand{\mg}{multigrid }
\journal{Journal of Computational Physics}

\begin{document}

\begin{frontmatter}


\title{Multigrid In Energy Preconditioner for Krylov Solvers}%\tnoteref{label1}}
\author{R. N. Slaybaugh\corref{cor1}}%\corref{cor1}
\ead{rns37@pitt.edu}
\author{T. M. Evans\fnref{label1}}
%\ead{evanstm@ornl.gov}
\author{G. G. Davidson\fnref{label1}}
%\ead{davidsongg@ornl.gov}
\author{P. P. H. Wilson\fnref{label2}}
%\ead{wilsonp@engr.wisc.edu}

\address{Department of Mechanical Engineering and Material Science, University of Pittsburgh, 605 Benedum Hall, 3700 O'Hara Street, Pittsburgh, PA 15261, USA\corref{cor1}}
\address[label1]{Radiation Transport Group, Oak Ridge National Laboratory, PO BOX 2008 MS6170, Oak Ridge TN 37831, USA\fnref{label1}}
\address[label2]{Department of Nuclear Engineering and Engineering Physics, University of Wisconsin - Madison, 419 ERB, 1500 Engineering Drive, Madison, WI 52706, USA\fnref{label1}}


\begin{abstract}
We have added a new multigrid in energy (MGE) preconditioner to the Denovo discrete-ordinates radiation transport code. This preconditioner takes advantage of a new multilevel parallel decomposition. A multigroup Krylov subspace iterative solver that is decomposed in energy as well as space-angle forms the backbone of the transport solves in Denovo. The space-angle-energy decomposition facilitates scaling to hundreds of thousands of cores. The multigrid in energy preconditioner scales well in the energy dimension and significantly reduces the number of Krylov iterations required for convergence. This preconditioner is well-suited for use with advanced eigenvalue solvers such as Rayleigh Quotient Iteration and Arnoldi.
\end{abstract}


\end{frontmatter}

% \linenumbers

\section{Introduction}
\label{sec:intro}
The steady-state Boltzmann equation for neutron transport finds ``where all the neutrons are'' in a nuclear system. The more accurately this is known, the more accurately new systems can be developed. Challenging transport problems today are three-dimensional with up to thousands $\times$ thousands $\times$ thousands of mesh points, use up to $\sim$150 energy groups, include accurate expansions of scattering terms, and are solved over many angular directions. Such discretizations create systems with $10^{9-10}$ coupled algebraic equations. High-fidelity, coupled, multi-physics calculations are the next generation of ``grand challenge'' problems for reactor analysis, requiring that the finely-resolved neutron flux be calculated quickly and accurately.

Very large computers, like the the Jaguar machine \cite{Jaguar2012}, are now available to perform such high-fidelity calculations. Access to such machines has changed the types of neutron transport problems that can be solved practically. However, most existing solution methods are not able to take full advantage of new computer architectures, or they have convergence properties that limit their usefulness for difficult problems.  

The deterministic transport equation is often solved with iterative methods, and Krylov methods are becoming an increasingly popular choice. Low-cost preconditioners that can reduce the number of iterations needed for convergence are invaluable. Benzi's 2002 survey paper on preconditioning techniques for large linear systems states ``it is widely recognized that preconditioning is the most critical ingredient in the development of efficient solvers for challenging problems in scientific computation''\cite{Benzi2002}. Ten years later, this is still true. The need for preconditioning is particularly high for Krylov methods because the memory required and cost per iteration increase dramatically with the number of iterations \cite{Benzi2002}.

This work aims to improve convergence behavior for 3-D transport calculations with a new preconditioner that can use leadership-class machines fully. A multigrid in energy preconditioner (MGE) that uses a Krylov solver to do the multigroup transport solves has been added to the code Denovo \cite{Evans2010}. The MGE preconditioner reduces the number of Krylov iterations in both fixed source and eigenvalue problems. This success has been demonstrated on large, realistic problems of interest. The preconditioner is particularly well-suited for the fine energy structures needed in challenging discrete-ordinates calculations. 

This paper is organized as follows: Section~\ref{sec:background} gives background information about solving the transport equation. Section~\ref{sec:precond} provides an overview of multigrid preconditioners used in the nuclear community. The MGE method, usage options, and how it is parallelized are presented in Section~\ref{sec:methods}. Finally, test results demonstrating the preconditioner's impact and the associated conclusions are presented in Sections \ref{sec:results} and \ref{sec:conclusions}, respectively.

%-----------------------------------------------------------------------------------------------------
%-----------------------------------------------------------------------------------------------------
\section{Background}
\label{sec:background}
To understand MGE and how it can enable the use of leadership-class computers, the mathematical details of the transport equation and how it is solved must be discussed first. 

%-----------------------------------------------------------------------------------------------------
\subsection{Transport Equation}
\label{sec:transport}  
The fixed source form of the steady state Boltzmann transport equation is
\begin{align}
   [\hat{\Omega} \cdot \nabla + \Macro(\vec{r}, E)] &\psi(\vec{r}, \hat{\Omega}, E)  = q_{ex}(\vec{r}, \hat{\Omega}, E)  \nonumber \\
   &+ \int_0^{\infty} dE' \int_{4\pi} d\hat{\Omega'} \:\Macro_{s}(\vec{r}, E' \to E, \hat{\Omega'} \cdot \hat{\Omega}) \psi(\vec{r}, \hat{\Omega'}, E') \:,
\label{eq:neutron transport}
\end{align}
\noindent where the quantities are at location $\vec{r}$, at energy E, and are traveling in direction $\hat{\Omega}$. The angular neutron flux in neutrons per unit length squared per steradian is $\psi(\vec{r}, \hat{\Omega}, E)$, which expresses the location of the neutrons in phase space; $q_{ex}$ is an external source; $\Macro_{x}$ is the cross section that represents the likelihood of interaction $x$ (here $s$ is for scattering, $f$ is for fission, and no subscript indicates all interactions) in units of inverse length \cite{Lewis1993}.

In the presence of fission, the external source is replaced by a fission source,
\begin{equation}
  \chi(E) \int_0^{\infty} dE' \:\nu \Macro_{f}(\vec{r}, E') \int_{4\pi} d\hat{\Omega'} \:\psi(\vec{r}, \hat{\Omega'}, E') \:.
  \label{eq:fission source}
\end{equation}
Here $\chi(E)$ is the spectrum specifying the energy distribution of neutrons born from fission, and $\nu$ is the average number of neutrons released per fission \cite{Lewis1993}. 

It is often of interest to know the asymptotic behavior of a system with fission. A reactor is called ``critical'' if the chain reaction is self-sustaining and time-independent. If the system is not in equilibrium, then the asymptotic neutron distribution, or fundamental mode, will grow or decay exponentially over time. A convenient way to capture this behavior is to assume $\nu$ can be adjusted to obtain a time-independent solution by replacing it with $\frac{\nu}{k}$, where $k$ is the parameter expressing the deviation from critical. A spectrum of eigenvalues can be found, but at long times only the non-negative solution corresponding to the largest real eigenvalue will dominate \cite{Lewis1993}. 

A brief aside about some properties of the transport equation will aid in understanding the solution techniques. In a void, the transport equation is like a hyperbolic wave equation. For highly-scattering regions where $\Macro_{s}$ is close to $\Macro$, the equation becomes elliptic for the steady-state case. If the scattering is forward-peaked then the equation is parabolic. All of these classes of equations have different solution strategies, so it is difficult to find one method that works well in all cases of physical interest \cite{Adams2002}.  

To numerically solve Equation~\eqref{eq:neutron transport} it is discretized in space, angle, and energy. This work uses the multigroup energy approximation, the scattering term is expanded in Legendre polynomials ($P_l$), and discrete ordinates (\Sn) are used to treat direction of neutron travel. There are many spatial differencing methods available, the discussion of which is beyond the scope of this document as the proposed work is not dependent upon the spatial discretization employed \cite{Evans2009d}. To ensure the new methods apply to the most general cases, it will be assumed that the matrices resulting from discretization are not necessarily symmetric. 

After all of the discretizations are performed, the multigroup \Sn equations can be written in operator form as
\begin{alignat}{2}
  \ve{L}\psi &= \ve{MS}\phi + q \:, \qquad &\text{(fixed source)} \label{eq:fxdsource} \\
  \ve{L}\psi &= \ve{MS}\phi + \frac{1}{k}\ve{M}\chi f^{T}\phi \:. \qquad &\text{(eigenvalue)} \label{eq:eigenvalue}
\end{alignat}
Here $\ve{L}$ is the first-order linear differential transport operator; $\ve{M}$ is the moment-to-discrete operator that projects the angular flux moments, $\phi$, onto discrete angles; $\ve{S}$ is the scattering matrix; $f$ contains the fission source, $\nu \Macro_{f}$; and $q$ is a source term. The angular flux moments are related to the angular flux through the discrete-to-moment operator: 
\begin{equation}
 \phi = \mathbf{D} \psi \:. \label{eq:disc to mom}
\end{equation}
Using this relationship, Equations \eqref{eq:fxdsource} and \eqref{eq:eigenvalue} can be rearranged such that they are a function of only $\phi$. The formulation is aided by defining $\ve{T} = \ve{DL}^{-1}$ and $\ve{F} = \chi f^{T}$ \cite{Evans2011}:
\begin{alignat}{2}
  (\ve{I} - \ve{TMS})\phi &= q \:, \qquad &\text{(fixed source)} \label{eq:OperatorFxdForm} \\
  (\ve{I} - \ve{TMS})\phi &= \frac{1}{k} \ve{TMF} \phi \:. &\text{(eigenvalue)} \label{eq:OperatorEvalForm}
\end{alignat}

Once the matrices are multiplied together, a series of single ``within-group'' equations that are each only a function of space and angle result. If the groups are coupled together by neutrons scattering from a low energy group to a higher energy group, then iterative ``multigroup'' solves over the coupled portion of the energy range may be required. If the eigenvalue is desired, an additional ``eigenvalue'' solve is needed \cite{Evans2009}. 

%-----------------------------------------------------------------------------------------------------
\subsection{Multigroup Solve}
\label{sec:multgroup}
Traditionally, the multigroup solve has been done with Gauss Seidel (GS). GS is iterative in energy. A space-angle solve using a within-group solver, such as source iteration or a Krylov method, is performed for each energy group in series. The groups are solved from $g=0$, the highest energy, to $g=G$, the lowest. For a group $g$ and an energy iteration index $j$ this is \cite{Evans2010}
\begin{equation}
  \bigl( \ve{I} - \ve{TMS}_{gg} \bigr) \phi^{j+1}_{g} = \ve{TM} \bigl( \sum_{g'=0}^{g-1}\ve{S}_{gg'}\phi^{j+1}_{g'} + \sum_{g'=g+1}^{G} \ve{S}_{gg'}\phi^{j}_{g'}  + q_{g} \bigr)  \:.
 \label{eq:up-GS}
\end{equation}

The first term on the right includes downscattering contributions from higher energies, and the second term represents upscattering contributions from lower energy groups that have not yet been converged for this energy iteration. Groups that only contain downscattering are simply solved once since the second term on the right is zero. Groups with upscattering, however, must be iterated until they converge. Convergence of GS is governed by the spectral radius of the system, so the method can be very slow when upscattering has a large influence on the solution \cite{Adams2002}. GS is fundamentally serial in energy because of how the group-to-group coupling is treated. 

A recently-added multigroup (MG) Krylov solver removes the traditional \\``within-group'' / ``multigroup'' iteration structure. This allows the solver to handle upscattering efficiently, and enables parallelization in the energy dimension. The solver has been shown to successfully scale to hundreds of thousands of cores. For example, a fixed-source test scaled from 69,102 cores to 190,080 cores with 98\% efficiency \cite{Slaybaugh2011}. 

The MG Krylov solver combines the space-angle and energy iterations to make one space-angle-energy iteration level. This allows the energy groups to be decomposed such that they can be solved in parallel. The space-angle-energy iterations are much like the within-group space-angle iterations, except that the iteration is over a block of groups instead of just one group. The block can include all groups or just upscattering groups, in which case the downscattering groups are treated in series as before. An added benefit of this solver is that Krylov methods generally converge more quickly than GS \cite{Trefethen1997}.

The multigroup Krylov method applied to the upscattering block is shown here, where $\ve{S}_{\text{up\_block}}$ contains the upscattering groups and $\ve{S}_{\text{up\_source}}$ has the downscattering only groups:
\begin{equation}
  \underbrace{(\ve{I} - \ve{TMS}_{\text{up\_block}})}_{\tilde{\ve{A}}}\phi_{\text{up\_block}}^{n+1} = \ve{TM}(\ve{S}_{\text{up\_source}}\phi_{\text{up\_source}}^{n+1} + q) \:.
  \label{eq:MGkrylov}
\end{equation}

Trilinos \cite{1089021} provides Denovo's Krylov solver, with a choice of either GMRES or BiCGSTAB \cite{Evans2010}. The Krylov solver is given an operator that implements the action of $\ve{\tilde{A}}$, or the matrix-vector multiply and sweep. In the MG Krylov solver, $\ve{\tilde{A}}$ is applied to an iteration vector, $v$, containing the entire upscattering block instead of just one group:
\begin{enumerate}
  \item matrix-vector multiply: $y = \ve{M}\ve{S}_{\text{up\_block}} v$,
  \item sweep: $z = \ve{T} y$,
  \item return: $v \leftarrow v - z$.
\end{enumerate}

To implement the energy parallelization, the problem is divided into energy sets, with groups distributed evenly among sets. After each set performs its part of the matrix-vector multiply, a global reduce-plus-scatter is the only required inter-set communication. Since each set uses the entire spatial mesh with the same spatial decomposition, the established performance of spatial scaling does not change. The space-angle decomposition in Denovo comes from the KBA wavefront algorithm \cite{Baker1998}. Evans et al.\ \cite{Evans2009d} showed scaling in space is limited to about 20,000 cores with a 500 million cell spatial mesh. 

The added energy decomposition offers the ability to further decompose a problem, even if the performance limit of spatial decomposition has been reached. The total number of cores is equal to the number of computational domains, that is, the product of the number of energy sets and the number of spatial blocks. For 20,000 spatial blocks and 10 energy sets, which is a reasonable decomposition, 200,000 cores will be used. See Ref. \cite{Evans2011} for more details. 

The MG Krylov solver makes it possible to use very large machines and fine energy structures. The addition of this solver is one of the key motivators for the MGE preconditioner because Krylov methods often need to be preconditioned to converge in a small number of iterations. A preconditioner that can scale well and handle many energy groups is essential for using the MG Krylov solver. 

%-----------------------------------------------------------------------------------------------------
\subsection{Eigenvalue Solve}
\label{sec:PI}
When the eigenvalue form of the transport equation is used, another iteration level is needed outside of the multigroup solves. A common way to solve the eigenvalue problem is with power iteration (PI). This method is attractive because it only requires matrix-vector products and two vectors of storage space. PI's convergence can be very slow for problems of interest, however, because the error is reduced in each iteration by a factor of $\ve{A}$'s dominance ratio, $\frac{\lambda_{2}}{\lambda_{1}}$, where the $\lambda$s are $\ve{A}$'s eigenvalues ordered as $|\lambda_{1}| > |\lambda_{2}| \ge \dots \ge |\lambda_{n}| \ge 0$.

For the standard form of power iteration, we define an energy-independent eigenvalue, $\Gamma = f^{T}\phi$. Then, Eq. \eqref{eq:OperatorEvalForm} can be written as
\begin{align}
  \ve{A}\Gamma &= k\Gamma \:, \qquad \text{where} \label{eq:EnergyIndepEval} \\
  \ve{A} &= f^{T}(\ve{I} - \ve{TMS})^{-1} \ve{TM}\chi \:. \nonumber
\end{align}
Note that this is the ordinary form of the eigenvalue problem with eigenpair $(\Gamma, k)$. The power iteration method, with $i$ as the iteration index, is then
\begin{equation}
  \Gamma^{i+1} = \frac{1}{k}\ve{A}\Gamma^{i} \:. \label{eq:PowerIteration}
\end{equation}
Inside of power iteration, the application of $\ve{A}$ to $\Gamma$ requires the solution of a multigroup problem that looks like a fixed source problem,
\begin{align}
  (\ve{I} - \ve{TMS})y^{i} &= \ve{TM}\chi \Gamma^{i} \:, \label{eq:EvalIndepFxdSource}\\
  \Gamma^{i+1} = f^T y^i \:. \label{eq:indepGamma}
\end{align}

Power iteration is just one of the eigenvalue solvers in Denovo. An Arnoldi solver and Rayleigh Quotient Iteration (RQI) are also available. For details about these methods and their implementation in Denovo see references \cite{Evans2011} and \cite{Slaybaugh2012}, respectively.

%-----------------------------------------------------------------------------------------------------
%-----------------------------------------------------------------------------------------------------
\section{Multigrid Preconditioning}
\label{sec:precond}
For many kinds of problems the total number of multigroup iterations required to adequately converge the transport equation is large. To reduce iteration count, a preconditioner can be applied to transform the system of interest into another equivalent system that has more favorable properties and is easier to solve. 
 
Right preconditioning leaves the right hand side of the equation unaffected and does not change the norm of the residual, which is used for convergence testing in most iterative methods. Right preconditioning is therefore often preferred over left or split preconditioning for iterative solvers \cite{Knoll2004}. A right preconditioner was implemented in this work, and the remaining discussion will be presented in right preconditioner format. 

Let $\ve{G}$ be a non-singular preconditioner. Then $\ve{A}x=b$ can be transformed as 
\begin{equation}
  \ve{AG}^{-1}y = b, \qquad  x = \ve{G}^{-1}y \:.
\end{equation}
The matrix $\ve{A}\ve{G}^{-1}$ is not formed in practice. The preconditioner can be applied by using some method to solve $\ve{G}y=c \to y \approx \ve{G}^{-1}c$, or by otherwise implementing the action of $\ve{G}^{-1}$ without explicitly forming and inverting $\ve{G}$ \cite{Benzi2002}.

There are many different types of preconditioners. These methods often are dependent upon the choice of spatial discretization employed, only apply to \\within-group iterations, or have other important limitations. As a result, preconditioning Krylov methods for solving the 3-D neutron transport problem is an active and vital area of research. We elected to investigate multigrid preconditioners with the goal of making a preconditioner that is independent of spatial discretization and can take advantage of the MG Krylov solver.

Multigrid methods are stationary iterative schemes that can be used alone or as preconditioners for other methods. Beginning in the late 1980s the nuclear community started using spatial \mg and/or angular \mg as both solvers and preconditioners. The first use of spatial \mg for transport equations in 1-D and 2-D was investigated by Nowak et al. Since that time, \mg has been used in multiple dimensions, for both isotropic and anisotropic scattering, and for various spatial discretizations \cite{Adams2002}. Some highlights from recent work are discussed below. All are applied to the \Sn neutron transport equation unless otherwise noted. 

In 1998, multigrid in space and multigrid in angle were used as preconditioners for Krylov methods and were tested for the 1-D, one-group, modified linear discontinuous (MLD) neutron transport equations by Oliveira and Deng~\cite{Oliveira1998}. They looked at isotropic scattering without absorption, isotropic scattering with absorption, and anisotropic cases. They had better results with multigrid than when using an Incomplete LU factorization method as a preconditioner.

In 2007, Chang et al.~\cite{Chang2007} used 2-D spatial \mg for the isotropic scattering case, with corner balance finite difference in space and a four-color block-Jacobi relaxation scheme. A bilinear interpolation operator and its transpose were used for grid transfer. The method had some trouble with heterogeneous problems. The authors assert their algorithm is parallelizable.

In 2010, Lee~\cite{Lee2010} developed a method to do \mg in space and angle simultaneously for two and three dimensions, isotropic and anisotropic scattering, one energy group, and a variety of spatial discretizations. The method can perform \mg in only space, only angle, or some combination thereof. It also handles thick and thin cells. 

\subsection{Two-Grid}
The two-grid acceleration method developed by Adams and Morel~\cite{Adams1993} was one of the earlier spatial \mg methods, and it has been built upon by others. It is intended to accelerate convergence of the outer iterations when upscattering is present. The outer iteration method is Gauss Seidel and acceleration is only applied to upscattering groups. The original work was done for slab geometries with a linear discontinuous (LD) discretization.

The two-grid method finds a correction for the isotropic component of the solution by solving a collapsed, one-group diffusion equation for the error. It is assumed that the zeroth moment of the error is a product of a spectral shape function and a space-dependent modulation function. To find the shape function, Fourier analysis is performed on GS. The zeroth moment of the cross sections is used to create the Fourier matrix, and then an eigenproblem is formed. Because the shape function is material-dependent, one such calculation is needed in each material region. Note that by using the zeroth moment, only the isotropic component of the solution is accelerated.

The diffusion equation in the two-grid approach differs slightly from the standard diffusion equation in that there is an extra term containing the gradient of the shape function. This is zero in homogeneous regions, but undefined at material interfaces. Adams and Morel found that neglecting the gradient term all together still gave good results for their test problems. This may not be true for more complex cases.

The two-grid method requires the diffusion operator to be consistent with the transport operator. This requirement can be difficult to meet for multi-dimensional problems, particularly for some spatial discretizations. 

\subsection{Two-Grid in Denovo}
\label{sec:TTG}
Denovo has a two-grid acceleration scheme based on the one developed by Adams and Morel. As mentioned, the iteration procedure uses a collapsed one-group diffusion equation to correct the low-order Fourier modes. Adams and Morel showed that the slowest converging spatial modes are diffusive and can be exactly computed in the infinite homogeneous case \cite{Adams1993}.

The original method can fail for systems of interest because of the consistency requirement for the discretization of the diffusion operator in multi-dimensional and multi-material problems. To make it applicable for the desired cases, the original two-grid method was modified by Evans et al.\ \cite{Evans2009d} by using a one-group transport equation instead of the diffusion equation. The modified method is called transport two-grid (TTG). The TTG method was derived to give the correct error estimation in the diffusive limit noted above. 

To execute the TTG scheme in Denovo, a transport sweep is conducted in each group, a residual is calculated, the low-order transport solve  is performed, an error form the transport sweep is calculated, and the scalar flux is updated. As with the original two-grid method, the TTG method is limited to correcting only the isotropic flux moments.

The additional cost of the method is like solving one extra group, so for cases with many upscattering groups the cost can be amortized. The acceleration equation can be preconditioned with Diffusion Synthetic Acceleration for additional speed. Finally, a reduced quadrature set can be used inside TTG to limit the cost. 

Some tests have shown TTG to be quite effective in improving the speed of convergence for upscatter problems when compared to unaccelerated GS \cite{Evans2009d}. While TTG solves the spatial discretization issue, it still only works with GS and the traditional within-group + multigroup iteration structure. 

%-----------------------------------------------------------------------------------------------------
%-----------------------------------------------------------------------------------------------------
\section{Multigrid In Energy}
\label{sec:methods}
Preconditioning is important for increasing the robustness of Krylov methods. This is particularly true for the multigroup Krylov solver. This solver can create large Krylov subspaces because it forms the subspaces with multiple-group-sized vectors. As a result, any reduction in iteration count will have a significant benefit in terms of memory and cost per iteration. Prior to this work, there was no preconditioner in Denovo that could work with the MG Krylov solver. 

The new preconditioner does multigrid in the energy dimension. To understand why multigrid in energy makes sense for neutron transport, some highlights about these methods are discussed here. See \emph{A Multigrid Tutorial} by Briggs, Henson, and McCormick~\cite{Briggs2000} for a thorough and approachable explanation of what multigrid methods are and why they work.

The error in $x_i$, the $i$th guess for $\ve{A}x_i=b_i$, can be written as a combination of Fourier modes. Each Fourier mode has a frequency, and the frequencies can range from low-frequency (smooth) to high-frequency (oscillatory). Iterative methods, also referred to as smoothers or relaxers, remove high-frequency error components quickly, but take many iterations to remove the low-frequency ones. 

The idea of multigrid methods is to take advantage of the smoothing effects of iterative methods by making smooth errors look oscillatory and thus easier to remove. Errors that are low-frequency on a fine grid can be mapped onto a coarser grid where they are high-frequency. A relaxer is applied on the coarser grid to remove the now oscillatory error components. The remaining error is mapped to a still coarser grid and smoothed again. The problem is restricted to coarser and coarser grids until it is on a grid which is coarse enough to directly invert the equations.

Next, the coarsest result is prolonged back to the next-finer grid and used to correct the solution there. A few relaxations are done on this finer grid. The errors are prolonged back up the chain, continuously correcting on finer grids, until the finest grid is reached. This entire process is called a V-cycle. 

Multigrid methods are differentiated by how many times the V-cycle is done and how many grids are used. The optimal combination of grids and cycles may depend on problem type. The addition of more grids and cycles will reduce error, but at an added cost.
 
A multigrid method was selected for many reasons. Multigrid methods remove the low-frequency error modes that require many Krylov iterations. The preconditioner was designed to take advantage of the energy decomposition used by the MG Krylov method. Each energy set can do work on its own grids and does not need to communicate with other energy sets. This is a communication savings compared to using grids in space or angle. An additional benefit is the simplicity of energy grids. Energy is one-dimensional, which allows for simpler coarsening and refinement than spatial or angular grids.

\subsection{Method}
\label{sec:mge}
Recall that right preconditioning is applied as $\ve{A} \ve{G}^{-1} \ve{G} \phi = b$, where \\$\ve{A} = \ve{I} - \ve{TMS}$. To implement this in Denovo, $y$ is defined as $\ve{G}\phi$ and the problem is broken into two steps: 
\begin{enumerate}
  \item with a Krylov method solve 
    \begin{equation}
      \ve{AG}^{-1}y = b \:, \label{eq:PrecondKrylov} 
    \end{equation}
  \item after finding $y$, calculate 
    \begin{equation}
      \phi = \ve{G}^{-1}y \:. \label{eq:PrecondPhi}
    \end{equation}
\end{enumerate}

To build energy grids, the energy group structure is coarsened so that each lower grid has fewer groups. The finest grid is the input energy structure, and the coarsest grid has one or a few groups. Each level has half as many groups as the previous level, rounded up if applicable. If there are $G+1$ groups on the fine grid there will be either $\frac{G+1}{2}$ or $\frac{G+2}{2}$ groups on the coarse grid. This is conceptually straightforward because the energy groups can be combined (restricted) and separated (prolonged) linearly. 
 
There are a variety of options that must be considered in designing a multigrid scheme: the restriction and prolongation operators, the relaxation method, the number and/or pattern of V-cycles to use, the number of relaxations to do on each grid level, and the depth of the V-cycle. Among these, the number of V-cycles, the number of relaxations per level, and the depth of the V-cycle have been implemented as user input options; the others are fixed. 

The implemented restriction operator is a simple averaging scheme. Neighboring fine data are averaged together to make coarse data.  
Recall that in multigrid methods the variables being restricted and prolonged are error modes. For a grid with spacing $h$ and a next-coarser grid $2h$, the errors are restricted as $e_{g}^{2h} = \frac{1}{2}(e_{2g}^{h} + e_{2g+1}^{h})$ for $g = 0,...,G_{end}$. 
If $G$ is even then $G_{end} = \frac{G}{2}$; if $G$ is odd then $G_{end} = \frac{G-1}{2}$ and the lowest energy group's datum is just copied. This scheme was chosen so that the thermal energy groups would be more resolved, which should improve accuracy for thermal reactors. The errors are restricted every time there is a transfer to a coarser grid.

The cross sections are restricted from the finest to the coarsest grid during problem initialization, and they do not change thereafter. The total and fission cross sections are restricted in the same way as the errors. Scattering is slightly more complicated since it has two indices, $g$ and $g'$. When there are an even number of groups all cross sections are treated the same way: for $g' = 0,..., G'_{end}$ and $g = 0, ..., G_{end}$,
\begin{align}
  \Sigma_s^{2h}(g,g') = &\frac{1}{4}[\Sigma_s^{h}(2g,2g') + \Sigma_s^{h}(2g+1,2g') \nonumber 
  \\ &+ \Sigma_s^{h}(2g,2g'+1) + \Sigma_s^{h}(2g+1,2g'+1)] \:. 
  \label{eq:XSSeven}
\end{align}
Unless there is upscattering in every group, some of the entries in Equation~\eqref{eq:XSSeven} will be zero. When there are an odd number of groups, the cross sections for the last $g$ and $g'$ groups are
  \begin{align}
    &\Sigma^{2h}_s(G_{end},G'_{end}) = \Sigma^{h}_s(2G,2G') \:, \\
    \text{for } &g' = 0,...,G'_{end} \:: \nonumber \\
    &\Sigma^{2h}_s(G_{end},g') = \frac{1}{2}[\Sigma^{h}_s(2G,2g') + \Sigma^{h}_s(2G,2g'+1)] \:,\\
    \text{for } &g  = 0,...,G_{end} \::& \nonumber \\
    &\Sigma^{2h}_s(g,G'_{end}) =  \frac{1}{2}[\Sigma^{h}_s(2g,2G') + \Sigma^{h}_s(2g+1,2G')] \:.
  \end{align}

The cross sections could be flux-weighted and re-restricted every time a new value for the flux were available, i.e.\ every new application of the preconditioner, to more accurately preserve physics. The experience of the computational community, however, has been that preconditioners do not need to rigorously and accurately preserve physics to be effective. The trade-off between the physics preservation and the performance of the preconditioner cannot be known \emph{a priori}. The assumption in this work is that the cost of recomputing the cross sections in every preconditioner application is greater than the benefit of more accurately representing the physics.

To prolong from a coarse to a fine grid, the points that line up between the grids are mapped directly: $e_{2g}^{h} = e_{g}^{2h}$ for $g = 0, ..., G_{end}$. To fill in the intermediate points on the fine grid, the adjacent coarse values are averaged: $e_{2g+1}^{h} = \frac{1}{2}(e_{g}^{2h} + e_{g+1}^{2h})$ for $g = 0, ..., G_{end}-1$. In the $e_{2g+1}^{h}$ case, $G_{end} = \frac{G-1}{2}$ when $G$ is even to preserve the direct transfer of points that line up between the grids. The errors are prolonged every time there is a transfer to a finer grid. Since cross sections do not change, they are never prolonged. 

There are other restriction and prolongation operators that are more rigorous and would preserve more accuracy when transferring between grids than those implemented. For example, a full weighting restriction operator would be more rigorous in combination with the current prolongation method than the current restriction operator is \cite{Briggs2000}. This change would be straightforward to implement. 

An example of a more accurate and complex prolongation operator is to use shape functions to prolong from a coarse grid to a fine grid. The shape functions could be based on the previous iterate of the fine-grid vector, some known desirable expansion, etc. Such a change would be more difficult to implement and require more research. 

These issues were not investigated in this work since the chosen operators, which were simple to code and check for correct implementation, were sufficient in practice. Any error added by the restriction and prolongation operators was largely removed by the relaxations. In addition, because this is a preconditioner, all of the pieces do not need to be rigorous. Complex grid transfer operators are unlikely to be worth the expense; however, it may be of value to investigate other simple operators. 

The user chooses the number of V-cycles done for each preconditioner application. One V-cycle proceeds from the finest grid to the coarsest grid and back to the finest. The input option specifies the number of V-cycles that are concatenated together, with a default of 2. Each additional V-cycle should remove more error, but has a computational cost. 

The depth of the V-cycle can also be specified by the user. The default behavior is determined by the number of groups, such that the grids will be coarsened until there is only one energy group. The number of grids needed is given by \cite{BinaryTree2012}
\begin{equation}
  \text{floor}\bigl( \log_{2}(G-1) \bigr) + 2 \:.
  \label{eq:NumGrids}
\end{equation}
How this is handled when using energy sets is discussed below. 

Some number of relaxations are performed on each level while traversing down and up the grids in a V-cycle. The number of relaxations per level is a user input choice with a default of 2. Performing more relaxations per grid should remove more error, but has a computational cost. The implemented relaxation method is weighted Richardson iteration. When applied to the transport equation, we obtain
\begin{equation}
  \phi^{m} = \bigr(\ve{I} + \omega(\ve{TMS} - \ve{I})\bigl)\phi^{m-1} + \omega b^{m-1} \:,
  \label{eq:relax}
 \end{equation}
where $\omega$ is a constant selected by the user that defaults to 1. 

An important principle is that the preconditioner is only attempting to approximately invert $\ve{A}$. It is therefore reasonable to use a less accurate angular discretization in the preconditioner than the rest of the code. For example, the whole problem may be solved at $S_{10}$, but the preconditioner could only use $S_{2}$. There is an input option to specify an angular quadrature set to use in the preconditioner different from the angular quadrature set used in the rest of the problem; the default is to use the same set in both. At this time, this option has only been implemented for vacuum boundary conditions. 

\subsection{Parallelization}
\label{sec:parallelization}
A key attribute of this preconditioner is that it is parallelizable in energy. This is because it uses the energy sets introduced by the MG Krylov solver. There are two ways to handle energy grids and energy sets together. One way is to restrict from $G$ groups down to $1$ group just as if there were no energy sets. This requires cross-set communication as soon as there are fewer groups than sets. This also causes some logistical difficulties related to what data are held by which sets at various points in the calculation. 

The other approach is to prohibit cross-set communication by having each set do its own ``mini'' V-cycle. Each set restricts, prolongs, and relaxes on only its own groups. This strategy requires there to be at least two groups on every set. With an unequal number of groups per set, there is a choice between forcing all sets to have the same grid depth or allowing those with more groups to have deeper Vs. The first option enforces energy load balancing between sets, while the second allows the sets with more groups to remove more error. 

The preconditioner is implemented so that all sets use the same grid depth. The benefit of load balancing is likely to be greater than having some sets use an extra grid. Thus, each set restricts to one or two group(s), giving approximately $num\_sets$ total groups across all sets at the coarsest level. The number of grids needed is determined by the set with the minimum number of groups, since it will be the first to reach a grid with one group. This modifies Equation~\eqref{eq:NumGrids} to be
\begin{align}
  num\_g_{min} &= \text{floor}\bigl(\frac{num\_groups}{num\_sets}\bigr) \:, \\
  num\_grids &= \text{floor}\bigl( \log_{2}(num\_g_{min}) \bigr) + 2 \:.
  \label{eq:multisetGrids}
\end{align}

The communication costs and logistical complications of the first strategy seem likely to overwhelm the benefit gained by going to one group instead of $num\_sets$ groups. The second strategy was chosen because it involves much less communication and overhead cost. The value of this choice will become apparent in the results section. 

With the implemented energy set strategy, there are trade-offs between the number of sets and the number of grids for a fixed number of groups. When there are more sets, more cores can be used at once and wall time should decrease. When there are fewer sets each V-cycle can go deeper, so the preconditioner should be more effective at reducing iteration count.  

\section{Results}
\label{sec:results}
Some tests were done to characterize the impact of preconditioning several problem types, as well as to investigate how different variables affect the preconditioner's behavior. The preconditioning parameters are the Richardson iteration weight, $w_{k}$, the number of V-cycles per preconditioner application, and the number of relaxations per level. 

The syntax used throughout this section will be that $w\#$ is the weight, $r\#$ is the number of relaxations per level, and $v\#$ is the number of V-cycles, e.g.\ $w1r1v1$ is one relaxation per level, one V-cycle, and a weight of 1. Using more preconditioning means using larger values of $w$ and/or $r$ and/or $v$. Other issues investigated are changing the depth of the V-cycle, using a different quadrature set inside the preconditioner than in the rest of the problem, and strong scaling. 

The goal of using the preconditioner is to improve convergence behavior of the multigroup Krylov solves. The best metric for measuring this is the total number of multigroup Krylov iterations used in a calculation because it is the most consistent and fair measure. The number of eigenvalue iterations is also compared for the eigenvalue tests. This is a point of interest rather than a measure of goal attainment. The total number of Krylov iterations is the best proxy for convergence behavior as it encompasses the work that is done within each eigenvalue iteration.  

Timing comparisons should be considered heuristically unless otherwise specified. In cases where the calculations were done on a single core, the machine was not dedicated to these calculations and times could vary if the same calculations were repeated. Some problems use an optimized version of the code and others use a debug version. These should give the same iteration count, but not necessarily the same relative times between problems. Further, little effort has been made to optimize the preconditioner for speed. Once the multigrid in energy preconditioner has been optimized for efficiency, the preconditioned times should decrease. How much improvement can be gained is a matter for future study. 

Unless otherwise noted, all test problems used a step characteristic spatial solver, level-symmetric angular quadrature, and the grid depth was determined using the default approach, Equation~\eqref{eq:NumGrids} or Equation~\eqref{eq:multisetGrids}. The Krylov solver was GMRES, which is set to limit the number of \mg iterations to 1,000 if the problem does not converge earlier. The convergence tolerances are noted for each problem. The tolerance for the multigroup solve is the convergence tolerance used by GMRES in Trilinos \cite{1089021}. The eigenvalue tolerance is used by PI to determine if the eigenvalue has converged. In Denovo, PI also checks the L2-norm and the infinity-norm of the difference in the fission source between iterations. The default L2-norm tolerance is 1.0 and infinity-norm tolerance is 0.01.

%-------------------------------------------------------------------------------------------------------
\subsection{Fixed Source Parameter Studies}
\label{sec:fxdsrcparam}
Some fixed source tests are discussed first, with the selection of preconditioning parameters informed by initial scoping tests that are not reported here. Fixed source calculations are particularly useful because the preconditioner can be studied apart from eigenvalue iterations. These tests used a debug version of Denovo on a single core.

The first test was a small, vacuum boundary problem. It contained one material and used ten groups, five upscattering groups, $P_{0}$ scattering, an $S_{4}$ level-symmetric angular quadrature set, a 3 $\times$ 3 $\times$ 3 grid, and a multigroup tolerance of 1 $\times$ 10$^{-6}$. The first three groups contained an isotropic source.

In one set of tests, the weight, $w$, was varied while the relaxations per level, $r$, and number of V-cycles, $v$, were both set to 1 ($w\#r1v1$). The results are in Table~\ref{table:FxdSrcTstVacWeight}. In this and all subsequent results an $r$, $v$, and/or $w$ of 0 corresponds to the unpreconditioned case. The table heading ``Krylov Iters'' indicates the number of Krylov iterations needed for convergence. 

It is clear that increasing the weight is initially beneficial, reducing the iteration count to 6 from the unpreconditioned 10. After a certain level, increasing weight begins to increase the number of iterations. The number of iterations needed to converge generally continued to increase above a weight of 1.9, and the problem failed to converge beginning with a weight of 2.1. Near a weight of 2 there was some behavior inconsistent with this trend. This issue is under investigation and is not material to the overall conclusions of this work. 

\begin{table}[!h]
\caption{Small Fixed Source Problem with Vacuum Boundaries, Weight Variation Study}
\begin{center}
\begin{tabular}{| c | c | c | c |}
\hline 
Weight & Relaxations & V-cycles & Krylov Iters \\[0.5ex]
\hline
0.0 & 0 & 0 & 10 \\
1.0 & 1 & 1 & 7 \\
1.1 & 1 & 1 & 6 \\
1.2 & 1 & 1 & 6 \\
1.3 & 1 & 1 & 6 \\
1.4 & 1 & 1 & 7 \\
1.5 & 1 & 1 & 8 \\
1.6 & 1 & 1 & 10 \\
1.7 & 1 & 1 & 13 \\
1.8 & 1 & 1 & 20 \\
1.9 & 1 & 1 & 265 \\
\hline 
\end{tabular} \\
\end{center}
\label{table:FxdSrcTstVacWeight}
\end{table}

Next, the number of relaxations per grid and the number of V-cycles were varied with the weight fixed at 1.0 ($w1r\#v\#$). These results are shown in Table~\ref{table:FxdSrcTstVacRV}. Recall that an $r$ of 3 means the preconditioner performed 3 relaxations per level, and a $v$ of 3 means the preconditioner conducted 3 V-cycles per application. 

\begin{table}[!h]
\caption{Small Fixed Source Problem with Vacuum Boundaries, Relaxation Count and V-cycle Variation Study}
\begin{center}
\begin{tabular}{|c| c| c| c|}
\hline
Weight & Relaxations & V-cycles & Krylov Iters \\[0.5ex]
\hline
0.0 & 0 & 0 & 10 \\
1.0 & 1 & 1 & 7 \\
1.0 & 2 & 2 & 4 \\
1.0 & 3 & 3 & 4 \\
1.0 & 4 & 4 & 4 \\
1.0 & 5 & 5 & 4 \\
1.0 & 6 & 6 & 4 \\
1.0 & 7 & 7 & 4 \\
1.0 & 10 & 10 & 4 \\
\hline 
\end{tabular}
\end{center}
\label{table:FxdSrcTstVacRV}
\end{table}

Initially, increasing $r$ and $v$ reduced the number of iterations needed for convergence. After enough preconditioning was done that only 4 iterations were needed, no additional amount of preconditioning reduced the iteration count further. 

Another calculation using a large amount of preconditioning, $w1.3r10v10$, also yielded 4 iterations. This both confirms that more preconditioning did not improve results, and demonstrates that a large amount of preconditioning does not cause breakdown. 

%-------------------------------------------------------------------------------------------------------
The previous problem configuration was modified with reflecting boundary conditions. Both weight and $r$/$v$ were varied. The multigroup tolerance and was again 1 $\times$ 10$^{-6}$. Results for varying the weight using 1 relaxation per level and 1 V-cycle are in Table~\ref{table:FxdSrcTstReflWeight}. The results for varying the number of relaxations per grid and the number of V-cycles using a weight of 1 can be seen in Table~\ref{table:FxdSrcTstReflRV}. A few $r$/$v$ variations using a weight of 1.3 instead of a weight of 1 are also in Table~\ref{table:FxdSrcTstReflRV}.

\begin{table}[!h]
\caption{Small Fixed Source Problem with Reflecting Boundaries, Weight Variation Study}
\begin{center}
\begin{tabular}{|c| c| c| c|}
\hline
Weight & Relaxations & V-cycles & Krylov Iters \\[0.5ex]
\hline
0.0 & 0 & 0 & 15 \\
1.0 & 1 & 1 & 11 \\
1.1 & 1 & 1 & 10 \\
1.2 & 1 & 1 & 10 \\
1.3 & 1 & 1 & 10 \\
1.4 & 1 & 1 & 11 \\
1.5 & 1 & 1 & 13 \\
1.6 & 1 & 1 & 19 \\
1.7 & 1 & 1 & 50 \\
1.8 & 1 & 1 & 1000 \\
\hline 
\end{tabular} \\
\end{center}
\label{table:FxdSrcTstReflWeight}
\end{table}

\begin{table}[!h]
\caption{Small Fixed Source Problem with Reflecting Boundaries, Relaxation Count and V-cycle Variation Study}
\begin{center}
\begin{tabular}{|c |c |c |c|}
\hline
Weight & Relaxations & V-cycles & Krylov Iters \\[0.5ex]
\hline
0.0 & 0 & 0 & 15 \\
1.0 & 1 & 1 & 11 \\
1.0 & 2 & 2 & 5 \\
1.0 & 3 & 3 & 3 \\
1.0 & 4 & 4 & 2 \\
1.0 & 5 & 5 & 2 \\
1.0 & 6 & 6 & 2 \\
1.0 & 7 & 7 & 2 \\
1.0 & 10 & 10 & 1 \\
\hline
1.3 & 2 & 2 & 4 \\
1.3 & 3 & 3 & 2 \\
1.3 & 4 & 4 & 2 \\
1.3 & 7 & 7 & 2 \\
\hline 
\end{tabular}
\end{center}
\label{table:FxdSrcTstReflRV}
\end{table}

The behavior of this problem was similar to the vacuum case. For this test a weight of 1.8 rather than 2.1 prevented convergence. Like the vacuum test, increasing the weight parameter above 1 decreased the number of iterations initially, and then the number of iterations increased with increasing weight. No anomolous behavior was observed for this problem. 

The $r$ and $v$ study showed two new things compared to the vacuum case. One is that with fewer iterations per level and fewer V-cycles, a weight of 1.3 was better than a weight of 1. The other was that the number of iterations could be decreased to 1 when a large amount of preconditioning was done. Whether and when 1 iteration can be reached for a given problem is related to how well Richardson iteration works for that problem's characteristics. 

These two sets of results show that increasing $r$ and/or $v$ decreases Krylov iteration count. They also indicate that using a small amount of weight can be beneficial, but a large amount is not. These two small and simple test cases provide a foundation for choosing parameters for testing with larger and more complex problems. 

%-------------------------------------------------------------------------------------------------------
\subsection{Fixed Source Angle Comparisons} 
For preconditioners to be effective, it is helpful if they are fast. One simple way to reduce preconditioning time is to use fewer angles for the preconditioner solves than for the overall problem. A half iron, half graphite, fixed source problem with an isotropic source assigned to all groups in all spatial cells was used to test this feature. The problem has vacuum boundaries, 27 energy groups, 13 of which have upscattering, and a $P_0$ scattering expansion. For this calculation a grid of 20 $\times$ 20 $\times$ 20 was used and the multigroup tolerance was 1 $\times$ 10$^{-6}$. The preconditioning parameters were $w1r2v2$. The overall problem was solved with an $S_{8}$ angular quadrature. 

To compare the effect of the reduced angular quadrature, the quadrature order within the preconditioner was varied from $S_{2}$ to $S_8$. The problem was spatially decomposed over four cores and an optimized version of Denovo was used. The results are shown in Table~\ref{table:ang}. 

\begin{table}[!h]
\caption{Iron Graphite Fixed Source Problem, Preconditioner Quadrature Order Study}
\begin{center}
\begin{tabular}{| c | c | c | c |}
\hline
Preconditioner \Sn & Time (s) & \% Diff Time & Krylov Iters \\[0.5ex]
\hline
8 & 221.66 & n/a & 17 \\
6 & 133.93 & 39.58 & 17 \\
4 & 98.67 & 55.49 & 17 \\
2 & 59.15 & 73.31 & 18 \\
\hline 
\end{tabular}\\
\end{center}
\label{table:ang}
\end{table}

For this problem, changing the number of solution directions inside the preconditioner had little or no impact on the number of iterations needed for convergence. Reducing the number of ordinates used inside the preconditioner reduced solution time in all cases, with $S_2$ giving the largest improvement: 73\% speedup. As expected, changing the preconditioner did not change the solution more than the convergence tolerance. Using fewer solve directions inside the preconditioner had a very positive impact: it reduced time but did not meaningfully change the solution.

%-------------------------------------------------------------------------------------------------------
\subsection{Multiple Energy Sets}
\label{sec:multisets}
Another important area of investigation was how the preconditioner faired when using multiple energy sets. This preconditioner was built to take advantage of multiple energy sets and be able to scale efficiently.

To investigate the effect of multigrid in energy with multiple energy sets on the Krylov iterations without worrying about the impact on an eigenvalue calculation, the iron-graphite fixed source problem was considered first. These data were calculated using an optimized version of the code. 

For this study, the spatial grid was increased to $50 \times 50 \times 50$ and an $S_{4}$ quadrature set was used instead of $S_{8}$. The unpreconditioned version was compared to one with $w1r2v2$ on 1 to 10 sets. Note that with 27 groups, 10 is the maximum number of energy sets compatible with the preconditioner. The problem was also decomposed over 2 $\times$ 2 spatial blocks. The calculations therefore used between 4 and 40 cores. 

The preconditioned calculation took 27 GMRES iterations while the unpreconditioned took 123 iterations, regardless of the number of energy sets used. Because the number of iterations did not change with sets, the only thing to compare is time. The focus is on relative change in time rather than absolute time since there is still room to optimize the preconditioner. 

\begin{figure}[!h]
    \begin{center}
      \includegraphics [width=0.83\textwidth, height=0.47\textheight] {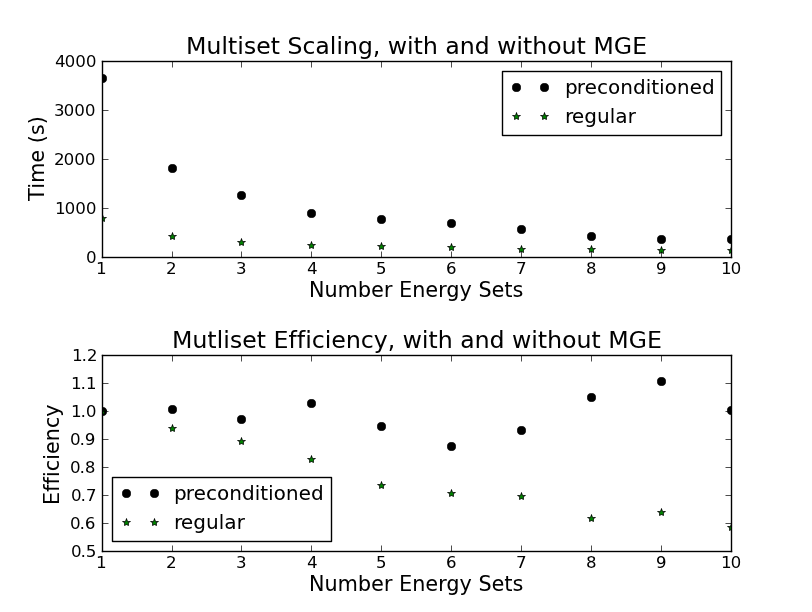}
   \end{center}
   \caption{Iron-Graphite Fixed Source Problem, Preconditioned Multiple Energy Set Study}
   \label{fig:FeC multisets}
\end{figure}

\begin{figure}[!h]
    \begin{center}
      \includegraphics [width=0.79\textwidth, height=0.46\textheight] {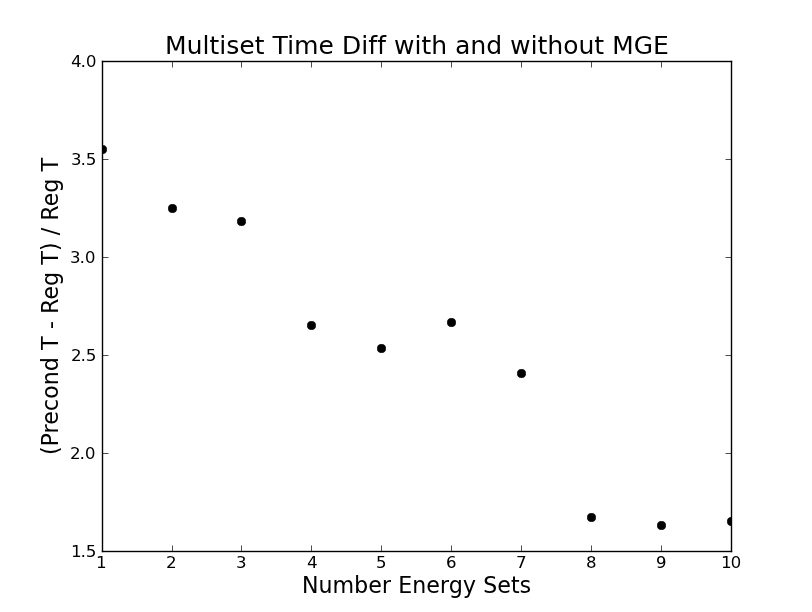}
   \end{center}
   \caption{Relative Time Difference In Iron-Graphite Fixed Source Multiple Energy Set Study}
   \label{fig:FeC multisets2}
\end{figure}

Two plots are shown in Figure~\ref{fig:FeC multisets}. The top plot shows the wall time for the preconditioned and unpreconditioned (``regular'') calculations as a function of number of energy sets. The bottom plot shows the efficiency of the regular and preconditioned tests, where the base case is 1 set. Figure~\ref{fig:FeC multisets2} is a plot of the relative difference between the two times. 

The 1-set wall time without preconditioning was $8.00 \times 10^{2}$ seconds and the 10-set time was $1.37 \times 10^{2}$ seconds. Another way to say this is that using 1 set took about 6 times longer than using 10 sets. If the problem scaled linearly in energy, it would have taken 10 times longer. The unpreconditioned efficiency degraded slowly with increasing energy sets, and the efficiency with 10 sets was less than 60\%.

With preconditioning, the 1-set time was $3.64 \times 10^{3}$ seconds and the 10-set time was $3.63 \times 10^{2}$ seconds. Going from 1 to 10 sets in this case was linear. The efficiency changed a bit with the number of energy sets, but ranged between about 90\% and 110\%. Thus, for some sets, the preconditioned case gave better than linear speedup. While the exact efficiencies may not be accurate because of timing variability, the general trend is clear. 

The difference between the efficiencies in the preconditioned and regular tests means that the preconditioned tests were accelerated more by using multiple energy sets than the unpreconditioned tests. As a result, the preconditioned times approached the regular times as more sets were used. This can be seen in the Figure~\ref{fig:FeC multisets2} plot of relative difference between the times. 

This test was modified and repeated for 1 set and 10 sets. The angular quadrature was $S_{8}$ and the problem was not decomposed in space. These choices were made to increase the runtime, ensuring that most of the reported time is solve time rather than setup or completion time. For the preconditioned cases, the preconditioning parameters were increased. Based on the first experiment, it is expected that the scaling will be the same or better when more preconditioning is used. 

With $w1.3r4v4$, an unnecessarily-large amount of preconditioning for this problem, but an amount good for investigative purposes, the number of Krylov iterations was 11. With 1 set this took $6.49 \times 10^{4}$ seconds to complete. With 10 sets this was reduced to $5.02 \times 10^{3}$ seconds. A 10-fold increase in computing power gave nearly a 13-fold decrease in run time, or an efficiency of 129\%. 

Without preconditioning, 124 Krylov iterations were required, the 1-set wall time was $8.81 \times 10^{3}$ seconds, and the 10-set time was $1.04 \times 10^{3}$ seconds. The ratio of 1-set to 10-set time was about 8.5, or 85\% efficient. This problem configuration scaled better than the previous configuration whether or not it was preconditioned. The comparison between the preconditioned and unpreconditioned cases show the same trends as the first energy scaling test. 

A similar scaling-in-energy test was done with an infinite medium, 27-group, eigenvalue problem. The same trends were observed. 
 
There are some clear reasons for the preconditioner's linear or super-linear scaling in energy. One is that this problem has a group structure that is not always balanced between sets. The multigrid preconditioner mitigates the penalty of energy-group load imbalance. Recall that each set uses the same number of grids, even if the number of groups per set is different. This causes the work in the preconditioner to be energy-load-balanced in all cases. Thus, the relative amount of time spent waiting because of load imbalance decreases when the preconditioner is used. 

The main reason for the preconditioner's good scaling is that as the number of sets increases, each application of the preconditioner becomes less costly. The total preconditioning cost goes down because the V-cycle becomes shallower. That is, each application of the preconditioner performs fewer total relaxations, and is therefore less time intensive. This effect becomes more pronounced with larger $r$ and $v$. If the total number of Krylov iterations remains constant with sets, then there is no trade-off, and the preconditioner simply costs less with increased energy parallelization.

An important outcome of this study is that the number of GMRES iterations did not change with the number of sets in the preconditioned cases. This means convergence improvement from the preconditioner does not come from the depth of the V-cycle, at least not for the problems tested. Only restricting down one or two grids had as much of an impact as restricting down something like six. This indicates it is not necessary to coarsen to one group. 

From an error mode reduction standpoint, this conclusion suggests what a Fourier expansion of the error in energy might look like. Because using a few somewhat-coarse grids had a large impact but using many much coarser grids did not, the bulk of the error might be intermediately oscillatory in energy. If the low-frequency error were dominant, then the coarsest grids would likely be necessary. Only solving on a fine grid, however, was not good enough, so the error is not only high-frequency either. This leaves the presence of modes that are between the two extremes as the likely culprit for slow convergence. 

The linear to super-linear energy scaling of the preconditioner warrants more discussion. Restricting down to and relaxing on the coarsest energy grids that were only created in the few-set cases did not provide convergence benefit. Much non-beneficial work was therefore done when only a few sets were used. When many sets were used and only a few grids that were not as coarse were created, all of this work was eliminated without any negative consequences. Thus, the energy scaling was very good.

These results show that limiting the number of energy grids to two or three is likely preferable for many problems since this strategy eliminates the non-beneficial work. By using fewer grids, the preconditioner will probably not cause the entire calculation to scale as well in energy. However, the preconditioner does not communicate between sets and it does load-balance in energy. The preconditioner should then, at worst, leave the energy scaling behavior unchanged. 

%-------------------------------------------------------------------------------------------------------
\subsection{V-cycle Depth}
\label{sec:vdepth}
The results of the multiple energy set study prompt an investigation of the V-cycle depth. The iron-graphite problem with a $50 \times 50 \times 50$ mesh was used for this purpose. The total problem quadrature was $S_8$, with $S_2$ used inside the preconditioner. The preconditioning parameters were $w1r2v2$. The problem was solved with an optimized version of Denovo on 25 cores on the orthanc cluster at ORNL, divided into 5 $x$-blocks, 5 $y$-blocks, 1 $z$-block, and 1 energy set. 

The depth of the V-cycle was varied from 1 to 6. A depth of 1 creates only 1 grid and a depth of 6 is the maximum possible for 27 energy groups. The test results are shown in Table~\ref{table:Vcycle}. The last row is the unpreconditioned case and the time difference is computed as $(t_{unPre}  - t)/t_{unPre}*100$. The minimum number of groups column indicates the minimum number of groups on the coarsest grid.
\begin{table}[!ht]
\caption{Iron Graphite Fixed Source Problem, V-cycle Depth Study}
\begin{center}
\begin{tabular}{| c | c | c | c | c |}
\hline
Depth & Min \# Groups & Solve Time (s) & \% Diff Time & Krylov Iters \\[0.5ex]
\hline
6 & 1 & 694.78 & 53.36& 28 \\
5 & 2 & 678.78 & 54.48 & 28 \\
4 & 4 & 643.52 & 56.80 & 28 \\
3 & 7 & 583.16 & 60.85 & 28 \\
2 & 14 & 462.77 & 68.94 & 28 \\
1 & 27 & 387.37 & 74.00 & 48 \\
n/a & n/a & 1489.69 & 0.00 & 124 \\
\hline 
\end{tabular}\\
\end{center}
\label{table:Vcycle}
\end{table}

The 1-grid case was the fastest. Increasing to 2 grids reduced the iteration count, though it also increased the time. Using more than 3 grids increased the solve time but did not change the iteration count. All cases took less time and used fewer Krylov iterations than the unpreconditioned case. 

This is the first problem presented in this paper that the preconditioned cases both used fewer iterations and took less time. Part of this success is that the preconditioner took less time per iteration because it used a reduced angular quadrature and a reduced grid depth. It was also the largest and most challenging of the previous test problems. 

These results confirm the multiple energy set study findings that using only a few grids is better than using many. The optimal number of grids will be problem dependent, since the shape of the error varies by problem. Without that knowledge \emph{a priori}, choosing a V-cycle depth of 2 is likely to give the best results. A V-cycle depth of 2 provides a significant reduction in iteration count while keeping the solve time low. 

%-------------------------------------------------------------------------------------------------------
%-------------------------------------------------------------------------------------------------------
\subsection{Eigenvalue}
\label{sec:eval}
All of the test problems presented have been fixed source. The MGE preconditioner will also be useful for eigenvalue calculations because it accelerates the multigroup solves in each eigenvalue iteration. Two eigenvalue test cases using Power Iteration are discussed here. These two tests are the most realistic so far. 

The eigenvalue tests have two main goals. One is to see if the preconditioner's general behavior is the same inside of an eigenvalue problem. The other is to see if the reduced Krylov iteration count is a beneficial trade-off for the cost of the preconditioner when used with Power Iteration. Because PI is such a simple solver, it may not be able to take advantage of the preconditioning as well as some more advanced eigenvalue solvers like Arnoldi or RQI.  

%-------------------------------------------------------------------------------------------------------
\subsubsection{3-D C5G7 Benchmark Study}
The preconditioner was applied to the 3-D version of the C5G7 mixed oxide (MOX) benchmark problem \cite{OECD-NEA2005}. This was solved with an optimized version of Denovo. The problem used three enrichment levels, an external mesh file, 10 materials, 7 groups, 4 upscattering groups, quadruple range angular quadrature using 4 polar and 6 azimuthal angles per octant, a $P_{0}$ scattering expansion, a multigroup convergence tolerance (or simply tolerance) of 1 $\times$ 10$^{-4}$, an eigenvalue tolerance (or $k$ tolerance) of 1 $\times$ 10$^{-5}$, and $k_{0}$ set to 1.14. The problem's dominance ratio is 0.7709, and the reference $k$ is 1.18655 $\pm$ 0.008. The problem was solved on the medium-sized OIC cluster at Oak Ridge, and each calculation was given 720 cores with 40 $x$-blocks, 18 $y$-blocks, and 5 $z$-blocks. The wall time limit was 12 hours. 

One drawback of this problem is that it is still fairly simple. The dominance ratio is not particularly high, so PI should not have much difficulty solving it. Further, there are not many energy groups, and the quadrature order is low. The primary goal of this calculation, though, was see whether the lessons learned about preconditioning parameters still hold in an eigenvalue problem. It is also good practice to test new methods on benchmark problems to verify they still calculate the correct solution.

The results are in Table~\ref{table:3-D c5g7}. The relative time is found by comparing the preconditioned wall time to the unpreconditioned wall time of $4.46 \times 10^{3}$ seconds. A number larger than 1.0 indicates the preconditioned case took longer. In all cases the computed $k$s were within the eigenvalue tolerance of one another, but none of them were within the uncertainty bounds of the reported benchmark; they were low by about 0.011. Subsequently, it was determined that using Gauss-Legendre product quarature with 6 polar and 8 azimuthal angles per octant with a linear discontinuous spatial solver instead of $QR$-$4$-$6$ and step characteristic gives the correct $k$. 
\begin{table}[!ht]
\caption{3-D C5G7 Benchmark, Preconditioning Parameter Scoping with Power Iteration}
\begin{center}
\begin{tabular}{| c | c | c | c | l | c |}
\hline
Weight & Relaxations & V-cycles & Krylov & PI & Rel Time$^{a}$ \\[0.5ex]
\hline
0.0 & 0 & 0 & 1,224  & 32 & 1.00 \\ %$4.46 \times 10^{3}$ \\
1.0 & 1 & 1 & 708    & 32 & 5.90 \\ %$2.12 \times 10^{4}$ \\
1.2 & 1 & 2 & 448    & 32 & 5.33 \\ %$2.38 \times 10^{4}$ \\
1.2 & 2 & 1 & 448    & 32 & 5.37 \\ %$2.39 \times 10^{4}$ \\
1.3 & 2 & 2 & 288    & 32 & 6.37 \\ %$2.84 \times 10^{4}$ \\
1.0 & 3 & 3 & 126    & 14$^{b}$  & 9.05 \\ %$4.04 \times 10^{4}$ \\
1.5 & 3 & 3 & 192    & 32 & 8.36 \\ %$3.73 \times 10^{4}$ \\
\hline 
\end{tabular}\\
$^{a}$compared to unpreconditioned, $4.46 \times 10^{3}$ seconds\\
$^{b}$tolerance = $1 \times 10^{-5}$, $k$ tolerance = $1 \times 10^{-3}$
\end{center}
\label{table:3-D c5g7}
\end{table}

The 3-D benchmark study shows the preconditioner with PI can reduce the total number of required Krylov iterations substantially. The number of iterations was cut nearly in half by using $w1r1v1$, and to one-third with $w1.2r1v2$. The $w1.2r1v2$ and $w1.2r2v1$ cases both took less time than $w1r1v1$. These two cases also took less time than the tests using more preconditioning. 

An intermediate amount of preconditioning gave the best trade-off between iteration count reduction and computational effort. Using less preconditioning cost less per iteration, but the extra iterations cost more time overall. Using more preconditioning reduced the iteration count, but the extra cost per iteration resulted in longer total times. 

The number of eigenvalue iterations for a given tolerance set were not changed by preconditioning. This means the preconditioner did not change the solution vector error in a way that had an impact on Power Iteration's behavior. 

These results continue to confirm that a small amount of weight works well. Increasing $r$ and $v$ decrease iteration count, but at what can be a high time penalty. The simplicity of the problem, the characteristics of using MGE with PI, or both, keep the preconditioner from being a time win for this test. The main conclusion is that an intermediate amount of preconditioning will likely provide the best balance of reduced iteration count for the time invested for problems that do benefit from the preconditioner. 

%-------------------------------------------------------------------------------------------------------
\subsubsection{Full PWR}
The next eigenvalue calculation was a whole-core pressurized water reactor (PWR)~\cite{Evans2011}. The core has a height of 4.2 m, assembly height of 3.6 m, and a lattice pitch of 1.26 cm. There are 157 fuel assemblies and 132 reflector assemblies. Each $17 \times 17$ assembly has $1/4$ symmetry and the fuel pins vary from 1.5\% to 3.25\% enrichment. The pins in the model are homogenized. 

In Denovo, each fuel pin has $2 \times 2 \times 100$ spatial cells. The calculations were done using an $S_{12}$ angular quadrature, $P_0$ scattering, and 44 energy groups, giving a total of 1.73 trillion unknowns. The PWR was decomposed over 137,984 cores on the Jaguar supercomputer using 578 $x$- and $y$-blocks, 700 $z$-blocks, and 11 energy sets. When the preconditioner was applied, the inside quadrature was $S_2$, the grid depth was 2, and $w1r2v2$ was selected. 

With tolerances of $1 \times 10^{-2}$, the unpreconditioned case converged in 24 eigenvalue iterations, using 901 Krylov iterations and taking 92.09 wall-minutes. When preconditioning was added, the eigenvalue iteration count was 28, the number of Krylov iterations was reduced to 312, but the time was increased to 250.09 wall-minutes. The preconditioner reduced the total number of required Krylov iterations by a factor of three, but took much longer. The Denovo-estimated dominance ratio for this tolerance is 0.938.

Both eigenvalue tests demonstrated that the multigrid in energy preconditioner dramatically reduces the Krylov iterations needed for convergence. Power iteration does not take advantage of this benefit. This motivates the need for advanced eigenvalue solvers.

Note that unpreconditioned PI converged these two reactor problems without trouble. We know that there are problems of interest that PI cannot solve quickly, or at all. These are the problems that require advanced eigenvalue solvers like RQI and Arnoldi, and will need preconditioners that can work with those solvers. The MGE preconditioner will likely be more useful for these solvers, and it is the only preconditioner available in Denovo that can operate on the multigroup Krylov solver.

A recent study of RQI in Denovo showed that it did not converge for challenging problems without preconditioning \cite{Slaybaugh2012}. This is a solver for which the MGE preconditioner could work well. For example, RQI could not solve the 3-D C5G7 benchmark without preconditioning. With preconditioning it behaved very well. There are other initial results showing the MGE preconditioner provides substantial improvement for RQI, and that preconditioned RQI converges more quickly than PI for difficult problems. Details about these results will be reported soon. 

%-------------------------------------------------------------------------------------------------------
%-------------------------------------------------------------------------------------------------------
\section{Conclusions}
\label{sec:conclusions}
With larger computers comes the potential to do more challenging neutron transport than ever before. To fulfill that potential, preconditioners that can scale efficiently are required. This paper discussed a new multigrid in energy preconditioner that can scale to hundreds of thousands of cores.

Preconditioning has been an important area of study in the transport community for many years. The use of a multigrid solver as a preconditioner for neutron transport began in the 1980s. Until now, multigrid methods have largely been conducted in the space and angle dimensions. This is the first multigrid in energy preconditioner for 3-D, large-scale, neutron transport. 

A variety of tests were done to characterize the MGE preconditioner implemented in Denovo. These tests showed that  performing 1 or 2 relaxations per V-cycle and using 1 or 2 V-cycles per preconditioner application provided the best overall combination of error reduction for the time invested. A weight of 1 when using the weighted Richardson method as the smoother was the most reliable choice. 

The preconditioner has the ability to use a lower-order quadrature than the rest of the problem, and this option should be selected whenever practical. Explicit investigation into the effect of limiting V-cycle depth showed that using a fairly shallow depth is sufficient. Truncating the V-cycle reduces calculation time while providing equivalent convergence benefit. The preconditioned cases converged more quickly and in fewer iterations than without preconditioning when these options were combined for the most difficult fixed source test.

The eigenvalue calculations with PI showed that the MGE preconditioner can significantly reduce the total number of Krylov iterations. For these two tests, however, using the preconditioner with PI was not faster than using PI alone. It may be that this preconditioner is not particularly well-suited to PI, or that it will only prove beneficial with more challenging problems. Initial results have shown that the preconditioner works quite well with RQI. 

The true potential value of the MGE preconditioner is its ability to capitalize on the new energy scaling capability. The multigroup Krylov solver has been shown to be extremely useful. It converges faster than Gauss Seidel, scales well in energy, and enables calculations with a new level of fidelity. MGE goes hand-in-hand with this solver.  

The new preconditioner can operate on entire multigroup vectors unlike other preconditioners in Denovo. MGE inherently provides load balancing across energy groups, and there is no inter-set communication during its application. Thus, it scales well in energy and can scale to hundreds of thousands of cores. It also facilitates the use of very fine energy structures. These are crucial traits for solving ``grand challenge'' transport problems. 

Overall, the MGE preconditioner reduces the number of Krylov iterations required. When parameters and options are chosen to reduce the time spent in the preconditioner, it can provide a time benefit as well. Important areas of future study are solving more difficult problems, using MGE with hundreds of energy groups, and studying its impact on advanced eigenvalue solvers like RQI.

\section{Acknowledgements}
\label{sec:acknowledements}
This research used resources of the Oak Ridge Leadership Computing Facility at the Oak Ridge National Laboratory, which is supported by the Office of Science of the U.S. Department of Energy under Contract No. DE-AC05-00OR22725. Additional thanks to the Rickover Fellowship Program in Nuclear Engineering sponsored by Naval Reactors Division of the U.S. Department of Energy. This fellowship sponsored the work from which this work is derived. 



%% References


\bibliographystyle{model1-num-names}
\bibliography{JCP_MGE}






% \end{thebibliography}

\end{document}